\documentclass[preprint,12pt]{elsarticle}

\usepackage{natbib}
\usepackage{color}
\usepackage{graphicx}
\usepackage{amssymb}
\journal{arXiv.org}

\begin{document}

\begin{frontmatter}

\title{Regional controllability of anomalous diffusion generated by the time fractional diffusion equations}

\author[dhu,ucmerced]{Fudong Ge}
\author[ucmerced]{YangQuan Chen\corref{cor}}
\ead{ychen53@ucmerced.edu}
\author[dhutwo]{Chunhai Kou}
\cortext[cor]{Corresponding author}

\address[dhu]{College of Information Science and Technology, Donghua University, Shanghai 201620, PR China}
\address[ucmerced]{Mechatronics, Embedded Systems and Automation Lab, University of California, Merced, CA 95343, USA}
\address[dhutwo]{Department of Applied Mathematics, Donghua University, Shanghai 201620, PR China}

\begin{abstract}
This paper is concerned with the investigation of the
regional controllability of the time fractional diffusion
equations. First, some preliminaries and definitions of regional
controllability of the system under consideration are introduced,
which promote the existence contributions on controllability analysis.
Then we analyze the regional controllability with minimum
energy of the time fractional diffusion equations on two cases:
$B\in \mathbf{L}\left(\mathbf{R}^m, L^2(\Omega) \right)$
and $B\notin \mathbf{L}\left(\mathbf{R}^m, L^2(\Omega)
\right)$. In the end, two applications are given to
illustrate our obtained results.
\end{abstract}
\par
\begin{keyword}
regional controllability; minimum energy; anomalous diffusion process;  fractional order
\end{keyword}
\end{frontmatter}

\section{Introduction}
\setcounter{section}{1}\setcounter{equation}{0}
Recently fractional diffusion equations (FDEs) have attracted increasing
 interests and a great deal of contributions have been given to both in
 time and space variables \cite{Podlubny,Kilbas,oldham,zhouyong,Rapaic}.
 And it is confirmed that the fractional approach to anomalous diffusion
 models is appealing  compared to other approaches. For instance, due to
the nonlocal and hereditary properties of fractional operators, the
anomalous diffusion models generated by FDEs are developed effectively
to describe transport process in complex dynamic system.

As we all know, the anomalous diffusion processes in real world are
essentially distributed and the continuous time random walk (CTRW),
governed by the waiting time probability density function (PDF) and
the jump length PDF, is a useful tool to describe this phenomenon
\cite{ralf,Friedrich,Meerschaert,Cartea}.
In addition, when the waiting time PDF and jump length PDF are
power-law and independent, the anomalous transport process can be
derived by the fractional partial differential equations, namely
fractional Fokker-Planck and Klein-Kramers equations \cite{ralf}.
And the time fractional diffusion equation models anomalous sub-diffusion
and its solutions are transition densities of a stable L$\acute{e}$vy
motion, representing the accumulation of power law jumps
\cite{Friedrich,Meerschaert,Cartea}.

Moreover, it can be easily seen that the control of anomalous diffusion
problem generated by FDEs can be reformulated as a problem  of analysis
of infinite-dimensional control system. However, in the case of diffusion
 systems, it should be pointed out that, in general, not all the states
can be reached. So in this paper, we first introduce some notations on the
regional controllability of  FDEs, i.e., the system under consideration is
 only exactly (or approximately ) controllable  on a subset of the state
 space,  which can be regarded as a generalization of integer-order diffusion systems
 \cite{AEIJAI,214,217}. Based on the semigroup theory \cite{pazy}, the regional
 controllability with minimum energy of time FDEs of two  different kinds of
 cases: $B\in \mathbf{L}\left(\mathbf{R}^m, L^2(\Omega) \right)$ and
 $B\notin \mathbf{L}\left(\mathbf{R}^m, L^2(\Omega) \right)$ are discussed.
 More precisely, when $B\in \mathbf{L}\left(\mathbf{R}^m, L^2(\Omega) \right),$
our main result is derived by utilizing the the Balder's theorem
\cite{Valadier}. And when $B\notin \mathbf{L}\left(\mathbf{R}^m, L^2(\Omega)
\right)$, the Hilbert Uniqueness Methods(HUMs), which were first introduced
 by Lions \cite{Lions}, are used to obtain the regional controllability
 with minimum energy of the system under consideration.

The rest of this paper is organized as follows: some concepts on regional
controllability are presented in the next section. In Section 3, our main
results on the regional controllability analysis of time FDEs
are given. Two  applications are worked out in the last section.
\section{Preliminaries}

Let $\Omega$ be an open bounded subset of $\mathbf{R}^n$ with smooth boundary $\partial\Omega$,
$Q=\Omega\times [0,T]$, $\Sigma=\partial\Omega\times [0,T]$. Let $L^p(0,T;\Omega)$ $(p\geq 1)$
be the space of Bochner integrable functions on $[0,T]$ with the norm $\|x\|_{L^p(\Omega)}=(
\int_0^T{\|x(s)\|_{\mathbf{R}^n}^p}ds)^{1/p}$  and consider the following abstract
fractional state-space system
\begin{eqnarray}\label{problem}
\begin{array}{l}
^C_0D^{\alpha}_{t}z(t)=Az(t)+Bu(t), ~~ z(0)=z_0\in D(A),
\end{array}
\end{eqnarray}
 where $t\in [0,T],$ $0<\alpha<1$, $z\in  L^2(0,T;\Omega)$, $^C_0D^{\alpha}_t$ is the
Caputo fractional derivative, $D(A)$ holds for the domain of the operator $A$  and
A generates a strongly continuous semigroup $\{\Phi(t)\}_{t\geq 0}$ on the Hilbert
state space $L^2(\Omega)$. In addition, $z_0\in L^2(\Omega)$, $u\in  L^2\left(0,T; \mathbf{R}^m\right)$
and $B: \mathbf{R}^m\to L^2(\Omega)$ is a linear operator to be
specified later.

Next, we will introduce some definitions and lemmas to be used in the sequel.

$\mathbf{Definition~~2.1.}$ \cite{Podlubny,Kilbas}\label{caputo} The Caputo fractional derivative of order
$\alpha>0$ of a function $z$ is given by
\begin{eqnarray}
^C_0D_t^\alpha z(t)=\left\{\begin{array}{l}
\frac {1}{\Gamma(n-\alpha)}\int_0^t(t-s)^{n-\alpha -1}\frac{\partial^n}{\partial s^n}z(s)ds,
\\\frac{\partial^n}{\partial t^n}z(t),~~~~~\alpha =n,
\end{array}\right.\end{eqnarray}
where $t\geq 0,$ $n-1< \alpha \leq n$, $n\in \mathbf{N},$ provided that the right side of $(2)$ is pointwise defined.

Let $\omega \subseteq \Omega $ be a given region of positive Lebesgue measure and
$z_T\in L^2(\omega)$$($the target function$)$ be a given element of the state
space.
\subsection{The case of $B\in \mathbf{L}\left(\mathbf{R}^m,L^2(\Omega)\right)$}
If $B\in \mathbf{L}\left(\mathbf{R}^m,L^2(\Omega)\right)$, i.e., $B$ is a bounded continuous operator
from $\mathbf{R}^m$ to $L^2(\Omega)$ and there exists a constant $M_B$ such that $\|B\|\leq M_B$.

  Based on the argument from the contribution \cite{zhouyong}, we get that

$\mathbf{Definition~~2.2.}$ \cite{zhouyong}\label{def2.3} For any given $u\in  L^2\left(0,T; \mathbf{R}^m\right),$
a function $z\in L^2\left(0,T; \Omega\right)$  is said to be a mild solution of the system
$(\ref{problem})$, denotes by $z(\cdot,u)$, if it satisfies
\begin{eqnarray}\label{transform}
z(t,u)=S_\alpha(t)z_0+\int_0^t(t-s)^{\alpha-1}{K_\alpha(t-s)Bu(s)}ds,
\end{eqnarray}
where
\begin{eqnarray}S_\alpha(t)=\int_0^\infty{\phi_\alpha(\theta)\Phi(t^\alpha\theta)}d\theta \end{eqnarray} and
\begin{eqnarray} K_\alpha(t)=\alpha\int_0^\infty{\theta\phi_\alpha(\theta)\Phi(t^\alpha\theta)}d\theta.\end{eqnarray}
Here $\{\Phi(t)\}_{t\geq 0}$ is the strongly continuous semigroup generated by $A$, $\phi_\alpha(\theta)=
\frac{1}{\alpha}\theta^{-1-\frac{1}{\alpha}}\psi_\alpha(\theta^{-\frac{1}{\alpha}})$ and  $\psi_\alpha$ is a probability density function defined by
\begin{eqnarray*}\psi_\alpha(\theta)=\frac{1}{\pi}\sum\limits_{n = 1}^\infty {( - 1)^{n - 1}\theta^{-\alpha n-1}\frac{\Gamma(n\alpha+1)}{n!}
\sin(n\pi \alpha)}, \theta \in (0,\infty)\end{eqnarray*}
 satisfying the following properties \cite{Mainardi}
 \begin{eqnarray}
\int_0^\infty{e^{-\lambda \theta}\psi_\alpha(\theta)}d\theta=e^{-\lambda^\alpha},~
\int_0^\infty{\psi_\alpha(\theta)}d\theta=1,~~\alpha\in (0,1)\end{eqnarray}
and
\begin{eqnarray}\label{2.4}
 \int_0^\infty{\theta^\nu
\phi_\alpha(\theta)}d\theta=\frac{\Gamma(1+\nu)}{\Gamma(1+\alpha\nu)}, ~\nu\geq 0.\end{eqnarray}

In order to prove our main results, the following hypotheses are needed.

$(S_1)$ The semigroup $\{\Phi(t)\}_{t\geq 0}$ generated by  operator A is uniformly bounded on
 $L^2(\Omega)$, i.e., there exists a constant $M>0$ such that \begin{eqnarray}\sup\limits_{t\geq 0}
 \|\Phi(t)\|\leq M.\end{eqnarray}

$(S_2)$ For any $t > 0,$ $\Phi (t)$ is a compact operator.

$\mathbf{Lemma~~2.1.}$ \cite{zhouyong}\label{2.2}

$(i)$ For any $t\geq 0$ the operator $S_\alpha(t)$ and $K_\alpha(t)$ are linear bounded operators, i.e.,
for any $x\in L^2(\Omega),$ \begin{eqnarray}\|S_\alpha(t)x\|_{L^2(\Omega)} \leq M\|x\|_{L^2(\Omega)}
\end{eqnarray} and \begin{eqnarray}\|K_\alpha(t)x\|_{L^2(\Omega)}\leq \frac{\alpha M}{\Gamma(1+\alpha)}
\|x\|_{L^2(\Omega)},\end{eqnarray}
where $M$ is defined in the inequality $(8).$

$(ii)$ Operators $\{S_\alpha(t)\}_{t\geq 0}$ and $\{K_\alpha(t)\}_{t\geq 0}$ are strongly continuous,
this is, for $\forall x\in L^2(\Omega)$ and $0 \leq t_1\leq t_2\leq T$, we have
\begin{eqnarray}\|S_\alpha(t_1)x-S_\alpha(t_2)x\|_{L^2(\Omega)}\to 0 \end{eqnarray} and
\begin{eqnarray} \|K_\alpha(t_1)x-K_\alpha(t_2)x\|_{L^2(\Omega)}\to 0 \mbox{ as } t_1\to t_2.\end{eqnarray}

$(iii)$ For $t>0,$ $S_\alpha(t)$ and
$K_\alpha(t)$ are all compact operators if $\Phi(t)$ is compact.

$\mathbf{Definition~~2.3.}$ \label{rcdefinition1}

$(a_1)$ The system $(\ref{problem})$ is said to be regionally exactly controllable(or $\omega-$exactly
controllable) if for any $z_T\in L^2(\omega)$, there exists a control $u\in L^2(0,T; \mathbf{R}^m)$  such that
\begin{eqnarray}
p_\omega z(T,u)=z_T.
\end{eqnarray}

 $(a_2)$ The system $(\ref{problem})$ is said to be regionally approximately  controllable(or
 $\omega-$approximately controllable) if for all $z_T\in L^2(\omega)$, given $\varepsilon>0,$
 there exists a control $u\in L^2(0,T; \mathbf{R}^m)$  such that
\begin{eqnarray}
\|p_\omega z(T,u)-z_T \|_{L^2(\omega)}\leq \varepsilon,
\end{eqnarray}
where $p_\omega: L^2(\Omega)\to L^2(\omega)$, defined by $p_\omega z=z|_\omega$, is a projection operator.

\subsection{ The case of $B\notin \mathbf{L}\left(\mathbf{R}^m, L^2(\Omega) \right)$}

If $B\notin \mathbf{L}\left(\mathbf{R}^m, L^2(\Omega) \right)$, similar to the argument in \cite{AEIJAI,214,217}, the extension definitions on regional controllability are introduced.

Take into  account that the system $(1)$ is line, without loss of generality, we suppose that $z_0=0$ in the following discussion.  Let $H:L^2(0,T;\mathbf{R}^m) \to  L^2(\Omega)$  be
 \begin{eqnarray}Hu=\int_0^T(T-s)^{\alpha-1}{K_\alpha(T-s)Bu(s)}ds,\forall u\in L^2(0,T;\mathbf{R}^m)_.\end{eqnarray}
  It follows from Definition 2.3 that the system $(1)$ is regionally approximately (exactly)
    controllable on $\omega$ if and only if
 \begin{eqnarray}\label{2.111}
\overline{im p_\omega H}=L^2(\omega)\left(\mbox{  respectively, }im p_\omega H=L^2(\omega)\right).
\end{eqnarray}
 $\mathbf{Definition~~2.4.}$

$(b_1)$ The system $(\ref{problem})$ is regionally exactly controllable if and only if
\begin{eqnarray}
ker p_\omega + im H= L^2(\Omega).
\end{eqnarray}

 $(b_2)$ The system $(\ref{problem})$ is said to be regionally approximately  controllable if and only if
\begin{eqnarray}\label{2.100}
ker p_\omega + \overline{im H}=L^2 (\Omega).
\end{eqnarray}

Suppose that $\{\Phi^*(t)\}_{t\geq 0}$, generated by the adjoint operator of $A$, is also a $C_0$
semigroup on the Hilbert state space $L^2(\Omega)$.
Then for any $v\in L^2(\Omega),$ it follows from $\left<Hu,v\right>=\left<u,H^*v\right>$ that
\begin{eqnarray}H^* v=B^*(T-s)^{\alpha-1}K_\alpha^*(T-s)v,\end{eqnarray} where $<\cdot,\cdot>$ is
the duality pairing of space $L^2(\omega)$, $B^*$ is the adjoint operator of $B$ and
$K^*_\alpha(t)= \alpha\int_0^\infty{\theta\phi_\alpha(\theta)\Phi^*(t^\alpha\theta)}d\theta.$ Then we have
$\overline{im p_\omega H}=L^2(\omega)$ is equivalent to
\begin{eqnarray}
ker H^* \cap im p_\omega^*=\{0\},
\end{eqnarray}
where $p^*_\omega:L^2 (\omega )\to L^2 (\Omega ),$
the adjoint operator of $p_\omega$, is
\begin{eqnarray}\label{pomega}
p_\omega^* z(x):=\left\{\begin{array}{l}z(x), ~~x\in \omega ,  \\0,~~~x\in \Omega\backslash \omega.\end{array}\right.
\end{eqnarray}
\section{Regional controllability analysis of the time FDEs}

In this section, we will explore the possibility of finding a
minimum energy control which steer the time FDEs $(\ref{problem})$ from the initial state $z_0$ to a target function $z_T$ on the region $\omega$.

Let $U_T=\{u\in  L^2\left(0,T; \mathbf{R}^m\right): p_\omega z(T,u)=z_T\}$.
Consider the following minimum problem
\begin{eqnarray}\label{minimum}
\inf\limits_u J(u)=\inf\limits_u \left\{{\int_0^T{\|u(t)\|^2_{\mathbf{R}^m}}dt}: u\in U_T\right\}. \end{eqnarray}
 \subsection{The case of $B\in \mathbf{L}\left(\mathbf{R}^m,L^2(\Omega)\right)$}

$\mathbf{Theorem~~3.1.}$\label{theorem3.1}
Suppose that $ B\in \mathbf{L}\left(\mathbf{R}^m,L^2(\Omega)\right)$ and the assumptions $(S_1)$,
$(S_2)$ hold, then the minimum problem $(\ref{minimum})$ admits at least one optimal solution
provided that the system $(\ref{problem})$ is $\omega-$approximately controllable.

$\textbf{Proof.}$
It is easy to see that $U_T$ is a closed and convex set. we first prove that the operator $H$ is
strongly continuous (see p.597, \cite{Hu}), which admits the existence of optimal control to the minimum problem $(22)$.
Moreover, according to the argument in \cite{Nieto},
since the operator $H$ is linear and continuous,  we only need to show that the operator is precompact.

For any $t\in [0,T], z_0\in L^2(\Omega)$, it follows from Lemma $\ref{2.2}$ that the term
$S_\alpha(t)z_0$ in Eq. $(\ref{transform})$ is strongly continuous. Let $N: L^2(\mathbf{R}^m)\to L^2(\Omega)$ be
\begin{eqnarray}\label{N}
Nu(t):=\int_0^t(t-s)^{\alpha-1}{K_\alpha(t-s)Bu(s)}ds,~~~t\in [0,T].
\end{eqnarray}
and we next show that $N$ is a relatively compact operator.

Let $\rho_r=\{ u\in  L^2(0,T; \mathbf{R}^m):\|u\|_{L^2(0,T; \mathbf{R}^m)} \leq r\}$.  For any fixed  $t\in(0,T]$, $\varepsilon,\delta\in(0,t)$,
$u\in \varrho_r,$  let
 \[\begin{array}{l}
 \widetilde{N}_{(\varepsilon,\delta)}u(t)=\alpha\int_0^{t-\varepsilon}\int_\delta^{\infty}(t-s)^{\alpha-1}\theta \phi(\theta)
 \Phi((t-s)^\alpha \theta)Bu(s)d\theta ds.\end{array}\]
 Since $\Phi(\varepsilon^q\delta)$ is compact and
\begin{eqnarray*}\widetilde{N}_{(\varepsilon,\delta)}u(t)
&=&\Phi(\varepsilon^\alpha\delta)\alpha\int_0^{t-\varepsilon}\int_\delta^{\infty}(t-s)^{\alpha-1}\theta \phi(\theta)
  \Phi((t-s)^\alpha \theta-\varepsilon^\alpha\delta)Bu(s)d\theta ds.\end{eqnarray*}
 we get that $\widetilde{N}_{(\varepsilon,\delta)}$ is relatively compact. Together with $\|Bu(\cdot)\|\leq M_Br<\infty,$
by $(i)$ in Lemma 2.1, for any $t\in [0,T]$, we have
\begin{eqnarray*}
\|\widetilde{N}u(t)-\widetilde{N}_{(\varepsilon,\delta)}u(t)\|
&=&\alpha\| \int_0^{t}\int_0^\delta{(t-s)^{\alpha-1}\theta \phi(\theta)
 \Phi((t-s)^\alpha \theta)Bu(s)}d\theta ds\\&~&
+\int_0^{t}\int_\delta^\infty{(t-s)^{\alpha-1}\theta \phi(\theta)
 \Phi((t-s)^\alpha \theta)Bu(s)}d\theta ds\\&~&-\int_0^{t-\varepsilon}
\int_\delta^\infty{(t-s)^{\alpha-1}\theta \phi(\theta)
 \Phi((t-s)^\alpha \theta)Bu(s)}d\theta ds\|\\
&\le&\alpha\left\| \int_0^{t}\int_0^\delta{(t-s)^{\alpha-1}\theta \phi(\theta)
 \Phi((t-s)^\alpha \theta)Bu(s)}d\theta ds\right\|\\&~&
+\alpha\left\|\int_{t-\varepsilon}^{t}\int_\delta^\infty{(t-s)^{\alpha-1}\theta \phi(\theta)
 \Phi((t-s)^\alpha \theta)Bu(s)}d\theta ds\right\|\\
&\le& MM_B r T^\alpha\int_0^\delta{\theta \phi(\theta)}d\theta+ \frac{MM_B r\varepsilon^\alpha }{\Gamma(1+\alpha)}
\to 0\end{eqnarray*} as $\varepsilon,
\delta\to 0$, where $M$ is defined in Eq. $(8).$ Then we conclude that $ N\varrho_r $ is a relatively compact set in $L^2(\Omega)$.

 Next, we shall prove that $Nu$ is equicontinuous  on $[0,T].$ For any $u\in \varrho_r,$ $0\leq\sigma_1<\sigma_2\leq T,$
 \begin{eqnarray*} &~&\|Nu(\sigma_2)-Nu(\sigma_1)\|\\
 &\leq&\left\|\int_0^{\sigma_1}[(\sigma_2-s)^{\alpha-1}-
 (\sigma_1-s)^{\alpha-1}]K_\alpha(\sigma_2-s)Bu(s)ds\right\|\\&~&+\left\|\int_0^{\sigma_1}(\sigma_1-s)^{\alpha-1}
 [K_\alpha(\sigma_2-s)-K_\alpha(\sigma_1-s)]Bu(s)ds\right\|\\&~&+\left\|\int_{\sigma_1}^{\sigma_2}(\sigma_2-s)
 ^{\alpha-1}K_\alpha(\sigma_2-s)Bu(s)ds\right\|\\ &\leq &\frac{MM_B r}{\Gamma(1+\alpha)}
 (\sigma_2^\alpha-\sigma_2^\alpha+(\sigma_2-\sigma_1)^\alpha)+A+\frac{MM_Br}{\Gamma(1+\alpha)}
(\sigma_2-\sigma_1)^\alpha\end{eqnarray*}
where $A=\left\|\int_0^{\sigma_1}(\sigma_1-s)^{\alpha-1} [K_\alpha(\sigma_2-s)-K_\alpha(\sigma_1-s)]Bu(s)ds\right\|$.
Since $\varepsilon>0$ small enough, we have
\begin{eqnarray*}
A&\leq &\int_0^{\sigma_1-\varepsilon}(\sigma_1-s)^{\alpha-1}\|K_\alpha(\sigma_2-s)-K_\alpha(\sigma_1-s)\|\|Bu(s)\|ds\\
&~&+\int_{\sigma_1-\varepsilon}^{\sigma_1}(\sigma_1-s)^{\alpha-1}\|K_\alpha(\sigma_2-s)-K_\alpha(\sigma_1-s)\|\|Bu(s)\|ds
\\&\leq&\left[\frac{M_Br}{\alpha}\left(\sigma_1^q
 -\varepsilon^q \right)\right]\mathop {\sup }\limits_{s\in[0,\sigma_1-\varepsilon]}\| K_\alpha(\sigma_2-s)-K_\alpha(\sigma_1-s)\|
 \\&~&+\frac{2MM_Br}{\Gamma(1+\alpha)}\varepsilon^{q}\\
 &\to& 0\end{eqnarray*}
 as $\sigma_2 \to \sigma_1$ due to the continuity of  $K_\alpha(t)(t>0)$ in the uniform operator topology.
It follows from the Arzela-Ascoli theorem \cite{Renardy} that the operator $N$ is precompact.
Thus, $H$ is strongly continuous, which guarantees the existence of optimal control
to the minimum problem $(\ref{minimum})$ under the fact that $U_T$ is a closed and convex set.

Further, if the system $(\ref{problem})$ is $\omega-$approximately controllable, for any $z_T\in \omega,$
suppose that $J(u^*)=\inf\limits_u J(u)=\varepsilon <\infty$, by the definition of infimum, we can
deduce that there exists a sequence $\{u_i\}_{i=1,2,\cdots}$ such that
$p_\omega z(T,u_i)=z_T,~u_i\in U_T\subseteq L^2\left(0,T; \mathbf{R}^m\right)$$(i=1,2,3,\cdots)$
and $J(u_i)\to J(u^*)$. Then we have $u_i\to_{_\omega} u^*$ in $L^2(0,T,\mathbf{R}^m).$

For any $t\in [0,T], $ by Definition 2.2 and Lemma 2.1, we get that
\begin{eqnarray*}&~~& \|p_\omega z(t,u^*)-p_\omega z(t,u_i)\|_{L^2(\Omega)}\\&= &\left\|p_\omega
\int_0^t(t-s)^{\alpha-1}{K_\alpha(t-s)B(u^*(s)-u_i(s))}ds\right\|\\
&\leq & \left\|\int_0^t{(t-s)^{\alpha-1}K_\alpha(t-s)B(u^*(s)-u_i(s))}ds\right\|
\\&\leq &  \frac{\alpha MM_B}{\Gamma(1+\alpha)}\int_0^t{(t-s)^{\alpha-1} \|u^*(s)-u_i(s))\|_{L^2(\mathbf{R}^m)}}ds,
\end{eqnarray*}
which yields that
\[p_\omega z(t,u_i)\to p_\omega z(t,u^*) \mbox{ in } C(0,T,\omega) \mbox { as } i\to \infty.\]
And since $U_T$ is closed and convex,  from Marzur Lemma \cite{Renardy} we see that $u^*\in U_T.$
Thus it follows from the Balder's theorem in \cite{Valadier} that
\[\varepsilon=J(u^*)=\lim \limits_{i\to \infty} J(u_i)\geq J(u^*)\geq \varepsilon,\]
which means that $u^* $ is the optimal solution of the minimum problem $(\ref{minimum})$. This completes the  proof.

\subsection{The case of $B\notin \mathbf{L}\left(\mathbf{R}^m,L^2(\Omega)\right)$}
If $B\notin \mathbf{L}\left(\mathbf{R}^m,L^2(\Omega)\right)$, for example, when the control
is pointwise or boundary control. The operator $N$ defined in Eq. $(23)$ is unbounded and then
$N$ is not relatively compact and  new methods should be introduced.

Here we will introduce the Hilbert uniqueness methods$($ HUMs $)$, which is first introduced
by Lions in \cite{Lions} to study the controllability problems of a linear distributed parameter
systems. Further, we note that this method is also available when $B$ is a bounded continuous operator.

Let $Z=im p_\omega H\subseteq  L^2(\omega)$, by duality $Z\subseteq  L^2(\omega)\subseteq Z^*$ and for any $f\in Z^*$, define
\begin{eqnarray}\|f\|_{Z^*}:=\int_0^T{\|B^*(T-s)^{\alpha-1}K_\alpha^*(T-s)p^*_\omega f\|^2}ds,
\end{eqnarray} where $p^*_\omega$ is defined in Eq. $(\ref{pomega})$.

$\mathbf{Lemma~~3.1.}$\label{lemma2.2}
 $\|\cdot\|_{Z^*}$ is a norm of space ${Z^*}$ provided that the system $(\ref{problem})$ is  $\omega-$approximately controllable.

 $\textbf{Proof.}$ If the system $(\ref{problem})$ is  $\omega-$approximately controllable,  we get that $\ker  H^* p^*_\omega=\{0\},$ i.e.,
\begin{eqnarray}B^*(T-s)^{\alpha-1}K_\alpha^*(T-s)p^*_\omega f=0 \Rightarrow f=0.\end{eqnarray}
 Hence, for any $f\in {Z^*}$, it follows from
  \begin{eqnarray*}
\|f\|_{Z^*}=\int_0^T{\|B^*(T-s)^{\alpha-1}K_\alpha^*(T-s)p^*_\omega f\|^2}ds=0
\Leftrightarrow B^*(T-s)^{\alpha-1}K_\alpha^*(T-s)p^*_\omega f=0\end{eqnarray*}
 that    $\|\cdot\|_{Z^*}$ is a norm of space $Z^*$ and the proof is complete.

Denote the completion of the set ${Z^*}$ with the norm $\|\cdot\|_{Z^*}$ again by $Z^*$.
 For each $f\in Z^*,$ since $f$ is a linear bounded functional on $Z,$ by the Riesz
 representation theorem, there exists a unique element in $Z$, denoted by $Pf$, such that
\begin{eqnarray}\label{3.17}
f(y)=(Pf,y)~~~~\mbox{ for all }y\in Z,
\end{eqnarray}
where $(,\cdot,)$ is the inner product in space $Z$.
Then we get that $P:Z^*\to Z$ is a linear operator and the following lemma holds.

$\mathbf{Lemma~~3.2.}$\label{lemma3.2}
The operator $P:Z^* \to Z$ is isometry.

\textbf{Proof.}
For any  $f\in Z^*,$  it follows from $(\ref{3.17})$ that
\[\|Pf\|_Z=\sup\limits_{\|y\|_Z=1}(Pf,y)=\sup\limits_{\|y\|_Z=1}\|f(y)\|=\|f\|_{Z^*}.\]
Then $\Re(P)\subseteq Z$ is  a closed subspace. To complete the proof, we should only
show that $\Re(P)= Z$. If not so, then there exists a $y_0\in Z,$ $y_0\neq 0$ such
that $(Px,y_0)=0.$ By $(\ref{3.17})$ , we have
\[f(y_0)=0~~~~\mbox{ for all }f\in Z^*,\]
which implies that $y_0=0,$ a contradiction.  Then we see that $\Re(P)= Z$ and the proof is complete.

Further, let $\wedge:Z^*\to Z$ be
\begin{eqnarray}\label{2.100}
\wedge f= p_\omega\varphi_1(T),
\end{eqnarray}
where $\varphi_1(t)$ is defined by
\begin{eqnarray}\label{2.11}
 \left\{\begin{array}{l}
^C_0D^{\alpha}_{t}\varphi_1(t)=A\varphi_1(t)+BB^*(T-t)^{\alpha-1}K_\alpha^*(T-t)f,
\\\varphi_1(0)=0 .
\end{array}\right.
\end{eqnarray}
Since for any $f\in Z^*,$ $y\in Z,$ by H$\ddot{o}$lder's inequality, we have
\begin{eqnarray*}
(\wedge f,y) &=&\int_{\Omega}p_\omega\int_0^T(T-s)^{\alpha-1}K_\alpha(T-s)B\times\\
&~&B^*(T-s)^{\alpha-1}K_\alpha^*(T-s)p_\omega^* f(x)dsy(x)dx\\
&\leq& \int_0^T{\|B^*(T-s)^{\alpha-1}K_\alpha^*(T-s)p^*_\omega f\|^2}ds\|y\|_Z\\
&\leq &\|f\|_{Z^*}\|y\|_Z
\end{eqnarray*}
and $\|\wedge f\|_Z\leq  \|f\|_{Z^*}$. Further, for any $f\in Z^*$, we have
\begin{eqnarray*}
(\wedge f,f)&=&\int_{\Omega}p_\omega\int_0^T(T-s)^{\alpha-1}K_\alpha(T-s)B\times\\&~&B^*(T-s)^{\alpha-1}K_\alpha^*(T-s)p_\omega^* f(x)dsf(x)dx\\
&=&\int_0^T\int_{\Omega}{\left[B^*(T-s)^{\alpha-1}K_\alpha^*(T-s)p_\omega^* f(x)\right]^2}dxds.
\end{eqnarray*}
Then if the system $(\ref{problem})$ is  $\omega-$approximately controllable on $[0,T]$, we get that $f=0.$ Thus it follows from the uniqueness of $P$ that $\wedge$ is an isomorphism from $Z^*$ to $Z$.

Next, suppose that $\varphi_0(t)$ satisfies
\begin{eqnarray}\label{2.11}
 \left\{\begin{array}{l}
^C_0D^{\alpha}_{t}\varphi_0(t)=A\varphi_0(t),
\\ \varphi_0(0)=z_0\in D(A),
\end{array}\right.
\end{eqnarray}
for all $z_T\in L^2(\omega),$  we have  $z_T=p_\omega\left[\varphi_1(T)+ \varphi_0(T)\right]$.
Further,
let $f$ be the solution of the following equation
\begin{eqnarray}\label{2.10}
\wedge f:= z_T- p_\omega\varphi_0(T).
\end{eqnarray}
Then we are ready to state the following theorem.

$\mathbf{Theorem~~3.2.}$\label{theorem3.2}
If the system $(\ref{problem})$ is $\omega-$approximately controllable, then for any
$z_T\in L^2(\omega),$ $(\ref{2.10})$ has a unique solution $f\in Z^*$ and the control
\[u^*=B^*(T-\cdot)^{\alpha-1}K_\alpha^*(T-\cdot)p_\omega^* f\] steers the system to $z_T$
 at time $T$ in $\omega$. Moreover, $u^*$ is the solution of the minimum problem $(\ref{minimum})$.

$\textbf{Proof.}$ By Lemma 3.1, we get that if the system $(\ref{problem})$ is
$\omega-$approximately controllable, then $\|\cdot\|_{Z^*}$ is a norm of space $Z^*$.
 Let the completion of $Z^*$ with respect to the norm $\|\cdot\|_{Z^*}$ again by $Z^*$.
  Then next we show that  the equation $(\ref{2.10})$ has a unique solution in $Z^*$.



 For any $f\in Z^*$,
by the definition of operator $\wedge$ in $(\ref{2.100})$, we get that
\begin{eqnarray*}<f,\wedge f>&=&<f, p_\omega\varphi_1(T)>\\&=&
<f,p_\omega\int_0^T(T-s)^{\alpha-1}{K_\alpha(T-s)Bu^*(s)}ds>\\&=&
\int_0^T{<f,p_\omega (T-s)^{\alpha-1}K_\alpha(T-s)Bu^*(s)>}ds \\&=&
\int_0^T{\|B^*(T-s)^{\alpha-1}K_\alpha^*(T-s)p_\omega^* f\|^2}ds \\&=&
\|f\|^2_{Z^*}.\end{eqnarray*}
Hence, it follows from Lemma $\ref{lemma3.2}$ and the Theorem 2.1 in \cite{Lions2}
that  the equation $(\ref{2.10})$ admits a unique solution in $Z^*$. Further,
 let $u=u^*$ in problem $(\ref{problem})$, we see that $p_\omega z(T,u^*)=z_T.$

For any $u_1,u_2\in L^2(0,T,\mathbf{R}^m )$ with $p_\omega z(T,u_1)=z_T$ and $p_\omega z(T,u_2)=z_T$,
we obtain that $p_\omega\left[z(T,u_1)-z(T,u_2)\right]=0.$  And for any $f\in Z^*,$ we have
\[<f,p_\omega\left[z(T,u_1)-z(T,u_2)\right]>=0.\]
It follows from
\[<p_\omega H(u_1-u_2), f>=<u_1-u_2,H^*p_\omega^*f>\]
that
\[\int_0^T{<u_1(s)-u_2(s),B^*(T-s)^{\alpha-1}K_\alpha^*(T-s)p_\omega^*f>}ds=0.\]
Then by
\begin{eqnarray*}J'(u_1)(u_1-u_2)&=&2\int_0^T{<u_1(s),u_1(s)-u_2(s)>}ds\\&=&2\int_0^T{<u^*(s),u_1(s)-u_2(s)>}ds\\&=&0,\end{eqnarray*}
we obtain that $u^*$ is the solution of the minimum problem $(\ref{minimum})$. This completes the proof.

\section{Example}
 In this section, we will introduce two examples which is reachable on $\omega$ but not reachable on the whole domain.

\textbf{Example 4.1.}
 Let us consider the following one dimension time FDEs with  $Bu=p_{[a_1,a_2]}u$,
 $0\leq a_1\leq a_2\leq 1$
\begin{eqnarray}\label{examplel}
\left\{\begin{array}{l}
^C_0D^{0.7}_{t}z(x,t)=\frac{\partial^2}{\partial x^2}z(x,t)+p_{[a_1,a_2]}u(t),~~~~[0,1]\times [0,T]
\\z(x,0)=z_0,{\kern 113pt}[0,1]
\\ z(0, t)= z(1, t)=0.~{\kern 80pt} [0,T]
\end{array}\right.
\end{eqnarray}
Corresponding to system $(\ref{problem})$, we have $A=\frac{\partial^2}{\partial x^2}$ and
\begin{eqnarray}
\Phi(t)z(x)=\sum\limits_{i=1}^{\infty}{\exp(\lambda_it)(z,\xi_i)_{L^2(0,1)}\xi_i(x)},~x\in [0,1],
\end{eqnarray}
where \[\lambda_i=-i^2\pi^2 ~\mbox{ and } ~\xi_i(x)=\sqrt{2}\sin (i\pi x),~x\in [0,1].\]
 Then  we get that the hypotheses $(S_1)$ and $(S_2)$ hold with $M=1.$
Further, we have
\begin{eqnarray}\begin{array}{l}
K_{0.7}(t)z(x)={0.7}\int_0^\infty{\theta\phi_{0.7}(\theta)\Phi(t^{0.7}\theta)z}d\theta\\
={0.7}\int_0^\infty{\theta\phi_{0.7}(\theta)\sum\limits_{i=1}^{\infty}{\exp(\lambda_it^{0.7}\theta)(z,\xi_i)_{L^2(0,1)}}\xi_i(x)}d\theta\\
={0.7}\sum\limits_{i=1}^{\infty}{(z,\xi_i)_{L^2(0,1)}\xi_i(x)\int_0^\infty{\theta\phi_{0.7}(\theta)\exp(\lambda_it^{0.7}\theta)}}d\theta.
\end{array}\end{eqnarray}
It follows from $(\ref{2.4})$ and the Taylor expansion of exponential function that
\begin{eqnarray}\begin{array}{l}
K_{0.7}(t)z(x)\\
={0.7}\sum\limits_{i=1}^{\infty}{(z,\xi_i)_{L^2(0,1)}\xi_i(x)\sum\limits_{j=0}^{\infty}{\int_0^\infty{\frac{(\lambda_it^{0.7})^j}{j!}
\theta^{j+1}\phi_{0.7}(\theta)}}}d\theta\\
=\sum\limits_{i=1}^{\infty}{(z,\xi_i)_{L^2(0,1)}\xi_i(x)\sum\limits_{j=0}^{\infty}{\frac{0.7(j+1)(\lambda_it^{0.7})^j}{\Gamma(1+{0.7} j+{0.7})}}}\\
=\sum\limits_{i=1}^{\infty}E_{0.7,0.7}(\lambda_it^{0.7})(z,\xi_i)_{L^2(0,1)}\xi_i(x),
\end{array}\end{eqnarray}
where  $E_{\alpha,\beta}(z):=\sum\limits_{i=0}^\infty
{\frac{z^i}{\Gamma(\alpha i+\beta)}},$ $~\mathbf{Re}{\kern 2pt}\alpha>0, ~\beta,z\in \mathbf{C}$
 is known as the generalized Mittag-Leffler function. Similarly, we have
\begin{eqnarray}
S_{0.7}(t)z(x)&=&\int_0^\infty{\phi_{0.7}(\theta)\Phi(t^{0.7}\theta)}d\theta\\
&=&\sum\limits_{i=1}^{\infty}{(z,\xi_i)_{L^2(0,1)}E_{0.7,1}(\lambda_it^{0.7})\xi_i(x)}.
\end{eqnarray}
What's more, since $A=\frac{\partial^2}{\partial x^2}$ generates a compact,
analytic, self-adjoint $C_0$ semigroup, we have
\begin{eqnarray*}\begin{array}{l}
(H^*z)(t)=B^*(T-t)^{-0.3}K_{0.7}^*(T-t)z(t)\\
=B^*(T-t)^{-0.3}\sum\limits_{i=1}^{\infty} E_{0.7,0.7}(\lambda_i(T-t)^{0.7})(z,\xi_i)_{L^2(0,1)}\xi_i(x)\\
=(T-t)^{-0.3}\sum\limits_{i=1}^{\infty} E_{0.7,0.7}(\lambda_i(T-t)^{0.7})(z,\xi_i)_{L^2(0,1)}\int_{a_1}^{a_2}{\xi_i(x)}dx.
\end{array}\end{eqnarray*}
Then it follows from $\int_{a_1}^{a_2}{\xi_i(x)}dx=\frac{\sqrt{2}}{i\pi}
\sin{\frac{i\pi (a_1+a_2)}{2}}\sin{\frac{i\pi (a_1-a_1)}{2}}$ that $ker H^*\neq
\{0\}$ $(\overline{im p_\omega H}\neq L^2(\omega))$ when $a_2-a_1\in Q,$ i.e.,
the system $(\ref{examplel})$ is not weakly controllable when $a_2-a_1\in Q$.
Thus, we can conclude that the system $(\ref{examplel})$ is not weakly controllable
on $[0,1]$ but on some appropriately subregion $[a_1,a_2]\subseteq [0,1]$ and
according to Theorem $\ref{theorem3.1}$, the minimum problem $(\ref{minimum})$ admits at least one optimal solution.

\textbf{Example 4.2.}
Consider the following time FDEs with pointwise control $Bu=u(t)\delta(x-b)$, $0< b< 1$, i.e., $B\notin \mathbf{L}\left(\mathbf{R}^m,L^2(\Omega)\right)$
\begin{eqnarray}\label{example2}
\left\{\begin{array}{l}
^C_0D^{0.7}_{t}z(x,t)=\frac{\partial^2}{\partial x^2}z(x,t)+u(t)\delta(x-b),~~~[0,1]\times [0,T]
\\ z(x,0)=0,{\kern 121pt}[0,1]
\\z(0, t)= z(1, t)=0.~{\kern 84pt} [0,T]
\end{array}\right.
\end{eqnarray}
 Here $Z=L^2(0,1),$ let $\omega=[\sigma_1,\sigma_2]\subseteq [0,1]$ and
if the system $(\ref{example2})$ is $\omega-$approximately controllable, since $A=\frac{\partial^2}{\partial x^2}$ generates a compact,
analytic, self-adjoint $C_0$ semigroup, similarly to the argument above, we have
\begin{eqnarray}\lambda_i=-i^2\pi^2 ,~ ~\xi_i(x)=\sqrt{2}\sin (i\pi x),,~x\in [0,1], \end{eqnarray}
\begin{eqnarray}
\Phi(t)z(x)=\sum\limits_{i=1}^{\infty}{\exp(\lambda_it)(z,\xi_i)_{L^2(0,1)}\xi_i(x)},~x\in [0,1]
\end{eqnarray}
and
\begin{eqnarray}
K_{0.7}(t)z(x)
=\sum\limits_{i=1}^{\infty}E_{0.7,0.7}(\lambda_it^{0.7})(z,\xi_i)_{L^2(0,1)}\xi_i(x),
\end{eqnarray}
Moreover, by Lemma 3.1, we get that if the system $(\ref{problem})$ is
$\omega-$approximately controllable,
\begin{eqnarray*}
&~&f\to \|f\|_{Z^*}\\
&=&\int_0^T{\left\|(T-s)^{-0.3}K_\alpha^*(T-s)p_\omega ^*f(b)\right\|^2}ds\\
&=&\int_0^T\left\|(T-t)^{-0.3} \sum\limits_{i=1}^{\infty}E_{0.7,0.7}(\lambda_i(T-t)^{0.7})(z,\xi_i)_{L^2(0,1)}
 p_\omega ^*f(b)\right\|^2ds\end{eqnarray*}
defines a norm on $Z^*$. It follows from Lemma 3.2 that
\begin{eqnarray}
\wedge f= p_\omega\varphi_1(T),
\end{eqnarray}
is a isometry form  $Z^*$ to $Z$,
where $\varphi_1(x,t)$ is the solution of the following equations
\begin{eqnarray}\label{2.11}
 \left\{\begin{array}{l}
^C_0D^{0.7}_{t}\varphi_1(x,t)=\frac{\partial^2}{\partial x^2}\varphi_1(x,)+(T-t)^{\alpha-1}K_\alpha^*(T-t)f(b),
\\ \varphi_1(x,0)=0.\\
 \varphi_1(0, t)=  \varphi_1(1, t)=0.
\end{array}\right.
\end{eqnarray}
Then by Theorem $\ref{theorem3.2}$, we see that
the control \[u^*(t)=(T-t)^{-0.3} \sum\limits_{i=1}^{\infty}E_{0.7,0.7}(\lambda_i(T-t)^{0.7})(z,\xi_i)_{L^2(0,1)}
 p_\omega ^*f(b)\] steers the system to $z_T$ at time $T$ in $\omega$, where $f$ is the solution of equations
\begin{eqnarray}
\wedge f=z_T- p_\omega\varphi_0(\cdot,T),
\end{eqnarray}
and $\varphi_0(t)$ solves
\begin{eqnarray}\label{2.11}
 \left\{\begin{array}{l}
^C_0D^{0.7}_{t}\varphi_0(x,t)=\frac{\partial^2}{\partial x^2}\varphi_0(x,t),
\\ \varphi_0(x,0)=z_0(x)\in D(A),\\
 \varphi_0(0, t)=  \varphi_0(1, t)=0.
\end{array}\right.
\end{eqnarray}
Moreover,  $u^*$ is the solution of the minimum problem $(\ref{minimum})$.

\section{CONCLUSIONS}

This paper is the first time to study the regional controllability analysis of the time fractional diffusion
equations on two cases: $B\in \mathbf{L}\left(\mathbf{R}^m, L^2(\Omega) \right)$
and $B\notin \mathbf{L}\left(\mathbf{R}^m, L^2(\Omega)
\right)$, which can be regarded as the extension of the existence contributions on controllability analysis of integer order
 \cite{AEIJAI,214,217}. The results we present here can also be extended to model real dynamic systems in complex dynamic system.
For instance, the problem of regional observability  of FDEs as well as the case of fractional super-diffusion equations with more complicated
dynamics are of great interest.\\

\noindent {\bf Acknowledgement.}
This work was supported  by Chinese Universities Scientific Fund No.CUSF-DH-D-2014061
 and the Natural Science Foundation of Shanghai (No. 15ZR1400800).


\bibliography{asme2e}

\begin{thebibliography}{12}
\bibitem{Podlubny} Podlubny I., 1999. "Fractional differential equations". {\it Academic Press, San Diego}.


\bibitem{Kilbas} Kilbas A. A., Srivastava H. M. and Trujillo J. J., 2006. "Theory and Applications of Fractional Differential Equations".
 {\it Elsevier}.

\bibitem{oldham}  Oldham K. B. and Spanier J., 1974. "The Fractional Calculus". {\it Academic Press, New York}.

\bibitem{zhouyong} Zhou Y. and Jiao F., 2010. "Existence of mild solutions for fractional neutral evolution equations".
{\it Comput. Math. Appl.}, 59, pp. 1063-1077.

\bibitem{Rapaic}  Rapai$\acute{c}$ M. R. and Jeli$\check{c}$i$\acute{c}$ Z. D., 2010. "Optimal control of a class of fractional heat diffusion systems".
{\it Nonlinear Dynam.}, 62, pp. 39-51.

\bibitem{ralf} Metzler R. and Klafter J., 2000. "The random walk's guide to anomalous diffusion: a fractional dynamics approach".
{\it Phys. Rep.}, 339, pp. 1-77.

 \bibitem{Friedrich} Friedrich R., Jenko F., Baule A. and Eule S., 2006. "Anomalous
diffusion of inertial, weakly damped particles". {\it Phys. Rev. Lett.}, 96, pp. 1-4.

\bibitem{Meerschaert} Meerschaert M. M., Nane E. and Vellaisamy P., 2013. "Transient anomalous sub-diffusion on bounded domains".
{\it Proc. Amer. Math. Soc.}, 141, pp. 699-710.

 \bibitem{Cartea} Cartea $\acute{A}$. and del-Castillo-Negrete D., 2007. "Fluid limit
 of the continuous-time random walk with general L$\acute{e}$vy jump distribution functions".
{\it Phys. Rev. }, E76, pp. 1-8.


\bibitem{AEIJAI} El Jai A. and Pritchard A. J., 1988. "Sensors and controls in the analysis of distributed systems".
{\it Halsted Press, New York}.

\bibitem{214}  El Jai A., Simon M. C., Zerrik E. and Pritchard A. J., 1995. "Regional controllability of distributed parameter systems".
{\it Internat. J. Control}, 62, pp. 1351-1365.

\bibitem{217} Zerrik E., Boutoulout A. and El Jai A., 2000. "Acutators and regional boundary controllability for parabolic systems".
{\it Internat. J. Systems Sci.}, 31, pp. 73-82.

\bibitem{pazy}  Pazy A., 1983.  "Semigroups of linear operators and applications to partial
 differential equations".  {\it Springer-Verlag, New York}.

\bibitem{Hu} Hu S. and Papageorgiou N.S., 1997. "Handbook of Multivalued Analysis (Theory)". {\it Kluwer Academic Publishers, Dordrecht Boston, London}.

\bibitem{Nieto} Debbouche A. and Nieto J. J., 2014.
"Sobolev type fractional abstract evolution equations with nonlocal conditions and optimal multi-controls".
{\it Appl. Math. Comput.}, 245, pp. 74-85.

\bibitem{Renardy} Michael R. and Robert C. R., 2004. "An introduction to partial differential
equations". Second edition. {\it Texts in Applied Mathematics, 13. Springer-Verlag, New York}.

\bibitem{Valadier} Valadier M., 1994. "Young measures, weak and strong convergence and the Visintin-Balder theorem".
  {\it Set-Valued Anal.}, 2, pp. 357-367.

\bibitem{Lions} Lions J.L., 1988. "Exact controllability, stabilization and perturbations for distributed systems".
{\it SIAM Rev.}, 30, pp. 1-68.

\bibitem{Mainardi} Mainardi F., Paradisi P. and Gorenflo R., 2000. "Probability distributions generated by fractional diffusion equations."
{\it in: J. Kertesz, I. Kondor (Eds.), Econophysics: An Emerging Science, Kluwer, Dordrecht}.

 \bibitem{Lions2} Lions J.L. and Stampacchia G., 1967. "Variational inequalities". {\it Comm. Pure Appl. Math.}, 20, pp. 493-519.

 \end{thebibliography}

\end{document}